\documentclass[twocolumn]{revtex4}
\usepackage{latexsym}
\usepackage{amsmath}
\usepackage{times}
\usepackage{amssymb}
\usepackage{fancyheadings}
\usepackage[T1]{fontenc}
\usepackage[utf8]{inputenc}
\usepackage{fancyhdr}
\usepackage{url}
\usepackage{hyperref}
\usepackage{color}
\usepackage{graphicx}

\begin{document}


\title{\bf A Kuramoto coupling of quasi-cycle oscillators}
\author{Priscilla E. Greenwood$^{1}$}
\author{Mark D. McDonnell$^{2}$}
\author{Lawrence M. Ward$^3$} \thanks{Corresponding author. Email: lward@psych.ubc.ca,  Tel.: +1 604 822 6309, Fax: +1 604 822 6923.}
\vspace*{3mm}
\affiliation{$^1$Department of Mathematics, University of British Columbia, Vancouver, BC, Canada}
\affiliation{$^2$Computational and Theoretical Neuroscience Laboratory, Institute for Telecommunications Research, University of South Australia, Mawson Lakes, SA, Australia}
\affiliation{$^3$Department of Psychology and Brain Research Centre, 2136 West Mall, University of British Columbia, Vancouver, BC, V6T 1Z4 Canada}
\date{\today}


\begin{abstract}
\noindent A family of stochastic processes has quasi-cycle oscillations if the oscillations are sustained  by noise. For such a family we define a Kuramoto-type coupling of both phase and amplitude processes. We find that synchronization, as measured by the phase-locking index, increases with coupling strength, and appears, for larger network sizes, to have a critical value, at which the network moves relatively abruptly from incoherence to complete synchonization as in Kuramoto couplings of fixed amplitude oscillators. We compare several aspects of the dynamics of unsynchronized and highly synchronized networks. Our motivation comes from synchronization in neural networks.

\vspace*{2ex}\noindent\textit{\bf Keywords}: Kuramoto coupling, quasi-cycles, neural oscillators, stochastic process, excitation-inhibition interaction, synchronization, Kuramoto model, phase locking index.
\\[3pt]
\end{abstract}

\maketitle

\thispagestyle{fancy}

\section{Introduction}

Strong bursty oscillations observed in, for example, EEG, or in local field potentials recorded from intracranial electrodes, derive from the activity of a number of subsystems of neurons that form an interactive network. In order to understand the dynamics of such a network, we would like to have a neural system model that produces similar synchronized oscillations. 

In~\cite{GMW14} we considered a stochastic neural system comprised of a group of pairs of excitatory and inhibitory neurons (E-I pairs). We described it in terms of a pair of interacting stochastic differential equations and found, using a new result about quasi-cycles, that such a model produces bursts of narrow-band oscillations, for example bursts of gamma-frequency oscillations.

The present paper is the first to study coupling of quasi-cycles. Our aim is to understand, for quasi-cycles, how the degree of synchronization depends on the strength of coupling of the subsystems.

Previous work has shown that groups of neurons with group potential oscillating independently can become synchronized when driven by a common input~\cite{Singer93,Varela01}. The paper~\cite{Wallace11} shows that both noisy limit cycles and quasi-cycles can arise from the Wilson-Cowan model, which is similar to the model we use to generate quasi-cycles. In~\cite{Daff11} a Kuramoto-type coupling is applied to deterministic Wilson-Cowan model limit cycles. In~\cite{Ton14} Kuramoto-type coupling is used with Freeman neural mass models that produce limit cycles with delay playing a role.

An aspect of our work that contrasts to~\cite{Wallace11, Daff11,Ton14} is that our stochastic subsystem models produce quasi-cycles, damped oscillations sustained by noise, and not limit cycles. The Wilson-Cowan model of interacting E (exciting) and I (inhibiting) populations includes a sigmoidal function whose role is to limit the primary variables in the dynamics to the interval [0,1]. A by-product is that the model has a limit cycle to which major attention has been directed. In~\cite{Daff11}, Kuramoto-type coupling is applied to a family of such deterministic systems, and the results refer to synchronization of limit cycles, in the presence of small equal amplitudes.

In contrast, here we omit the sigmoidal function in the initial model so that its centered dynamics for the model without noise would be simply oscillations damped to a fixed point at (0,0). We use a recent stochastic dynamics result of Baxendale and Greenwood~\cite{Bax11} to move to an equivalent model where amplitude and phase of quasi-cycles can be conveniently viewed as the primary variables satisfying a system of stochastic differential equations (SDEs). Then we introduce Kuramoto-type coupling to the amplitude and phase equations separately. The resulting model tells us how the degree of coupling affects the degree of synchronization among quasi-cycles.

Our results apply in the following general context: consider a system of stochastic oscillators, where each subsystem is described by a pair of linear, or locally linear, SDEs. Suppose the deterministic version of this system, with zero noise, has oscillations damped to a fixed point. Then, each stochastic subsystem has sustained oscillations. We coupled these stochastic subsystems and studied their synchronization. 

\subsection{Background I: A Class of Linear Noise Models called `quasi-cycle oscillators'}

We first write a linear model similar to that of~\cite{Kang10}, which in neuroscience is often called a Wilson-Cowan-type model. The Excitatory-Inhibitory (E-I) model for each subset of neurons (E-I pair or stochastic oscillator) appears as
\begin{align}\label{Kang1}
\begin{split}
& \tau_E dV_{E}(t) =(-V_{E}(t)+S_{EE}V_{E}(t)-S_{EI}V_{I}(t)) dt\\
& + \sigma_{E} dW_{E}(t)\\
& \tau_I dV_{I}(t)=(-V_{I}(t)-S_{II}V_{I}(t)+S_{IE}V_{E}(t)) dt\\
& + \sigma_{I} dW_{I}(t).
\end{split}
\end{align}
Here $W_{E}, W_{I}$ are independent, standard Brownian motions. The parameters $S_{EE}, S_{II}, S_{IE}, S_{EI}\ge 0$, are constants that represent the mean efficacies of the excitatory or inhibitory synaptic connections to post-synaptic neurons within each separate population, as indicated by the notation, with $S_{IE}$ representing input to inhibitory from excitatory neurons. These parameters, along with the time constants, $\tau_E, \tau_I$, and amplitudes of the Brownian motions, $\sigma_E, \sigma_I$, determine the oscillatory behaviour of the system and in particular its resonant frequency of oscillation. We limit our discussion here to the parameter ranges where the oscillation is narrow-band and thus has a distinct phase even though it arises from a stochastic process~\cite{GMW14}. Note that~\eqref{Kang1} can be interpreted as applying to a single pair of neurons, or even to the subthreshold dynamics of a single stochastic neuron~\cite{Engel08,Dit12}. Here, however, we consider that as a single oscillatory system characterized by a particular resonant frequency, regardless of how many neurons are involved in creating that system.  Note also that we can consider~\eqref{Kang1} to be a linear approximation, in the neighborhood of a fixed point, to a more elaborate dynamical system. 

System~\eqref{Kang1} has already been centered at the fixed point, (0,0), and can be written as what recently has been termed a linear noise model (or linear noise approximation if the original model is not linear),
\begin{equation}\label{Kang3}
d\mathbb{V}=-\mathbb{AV} dt + \mathbb{N} d\mathbb{W}
\end{equation}
where $\mathbb{V}=(V_E(t),V_I(t))^\top$, $d\mathbb{W}=(dW_E(t),dW_I(t))^ \top $,
\begin{equation}\label{V_vector}
\mathbb{A}=\begin{pmatrix}
			(1-S_{EE})/\tau_E& S_{EI}/\tau_E\\
			-S_{IE}/\tau_I& (1+S_{II})/\tau_I
			\end{pmatrix},\quad
\end{equation}
and
\begin{equation}
\mathbb{N}=\begin{pmatrix}
			\sigma_E/\tau_E & 0\\
			0 & \sigma_I/\tau_I
			\end{pmatrix}.
\end{equation}
			
When we take $\sigma_E=0$ and $\sigma_I=0$, the deterministic system obtained from Eqn.~(\ref{Kang3}) has oscillations that damp to the stable point (0,0) at a rate $\lambda$. In fact, Cowan's early work~\cite{Cowan70} used the fact that asymmetric coupling of neurons would typically generate oscillations because of the  presence of complex eigenvalues near the steady state. Under certain conditions the stochastic system in~(\ref{Kang3}), with elements of $\mathbb{N}$ nonzero, has sustained oscillations, called quasi-cycles, of a narrowband nature~\cite{Bressloff10,Bax11}. This is the case when the eigenvalues of the matrix $-\mathbb{A}$ are complex, $-\lambda \pm i \omega_d$, with $0<\lambda\ll\omega_d$. Here, as in~\cite{Kang10}, $\lambda$ is the damping rate of the oscillation,
\begin{equation}\label{Damping_rate}
\lambda=0.5\left[\frac{1-S_{EE}}{\tau_E}+\frac{1+S_{II}}{\tau_I}\right],
\end{equation}
and $\omega_d$ is its natural frequency,
\begin{equation}\label{omega_d}
\omega_d=\sqrt{\frac{S_{EI}S_{IE}}{\tau_E\tau_I} - 0.25\left[\frac{1-S_{EE}}{\tau_E}-\frac{1+S_{II}}{\tau_I}\right]^2},
\end{equation}
which is positive when the eigenvalues are complex. We give the name \emph{quasi-cycle oscillator} to a linear noise model~\eqref{Kang3} such that $\lambda$ and $\omega_d$ satisfy $0<\lambda<<\omega_d$, and hence has oscillations sustained by noise.

\subsection{Background II: An Asymptotically Equivalent Amplitude and Phase Model}
When a system~\eqref{Kang3} is a quasi-cycle oscillator, i.e., the damping rate is substantially less than the frequency of the oscillations, and the eigenvalues of $-\mathbb{A}$ are complex with negative real part, the damped oscillations of the deterministic version of the system will be sustained by `small' noise. In a neural system, this noise is always present because of fluctuating inputs occurring at synapses and because of the stochastic effects of ion channels~\cite{Schneidman.98}. For reasonable parameter values in \eqref{Kang1}, we have $0<\lambda\ll\omega_d$, i.e. the range of parameter values is such that  sustained oscillations are produced. This choice is justified in, e.g., \cite{GMW14} and \cite{Kang10}. Under these conditions, and according to a limit theorem of Baxendale and Greenwood \cite{Bax11}, the process (\ref{Kang3}) is approximated as $\lambda/\omega_d$ approaches 0, by 
\begin{equation}\label{Rotation}
\mathbb{V}(t)=[V_E(t), V_I(t)]^\top \approx \mathbb{V}^*(t) := \frac{\sigma}{\sqrt{\lambda}} \mathbb{QR}_{-\omega_d t} \mathbb{S}(\lambda t),
\end{equation}
 where $\mathbb{R}$ is the rotation
\begin{equation}\label{rota}
\mathbb{R}_s = \begin{pmatrix}
				\cos{(s)}&  -\sin{(s)}\\
				\sin{(s)}&   \cos{(s)}
				\end{pmatrix},
\end{equation}
$\mathbb{Q}$ is a matrix that transforms $\mathbb{A}$ in Eqn.~(\ref{Kang3}) into a canonical form
\begin{equation}\label{Q}
\mathbb{Q}^{-1}(-\mathbb{A})\mathbb{Q} = \begin{pmatrix}
									-\lambda&  \omega_d\\
									-\omega_d&  -\lambda
									\end{pmatrix}
									:= \mathbb{A}_1,
\end{equation}
$\mathbb{S}(t)$ is a standard two-dimensional Ornstein-Uhlenbeck process with independent components, and 
\begin{equation}\label{sigma}
\sigma = \sqrt{0.5\mathrm{Tr}(\mathbb{Q}^{-1}\mathbb{N}\mathbb{N}^\top(\mathbb{Q}^{-1})^\top)}
\end{equation}
is a scalar.

The approximation we used to obtain~\eqref{Rotation} from the linear model~\eqref{Kang3} is very different from the familiar method involving expansion and dropping of higher order terms. Thus we here sketch a proof of the result leading to~\eqref{Rotation} starting from the linear noise approximation~\eqref{Kang3}, keeping in mind that $0 < \lambda \ll \omega_d$ will be used in an essential way~\cite{Bax11}.

We make three changes of variables in~\eqref{Kang3}. First we transform the matrix $-\mathbb{A}$ to normal form in order to see more clearly the separate effects of the relatively slow damping characterized by $\lambda$, and the relatively fast rotation characterized by $\omega$. Using a matrix $\mathbb{Q}$ as in~\eqref{Q} (see~Eqn.~(4.4) in \cite{GMW14}), we write $\mathbb{Y}(t)=\mathbb{Q}^{-1} \mathbb{V}(t)$. Then
\begin{align}\label{dY}
d\mathbb{Y}(t)=\mathbb{A}\mathbb{Y}(t)-\mathbb{C}d\mathbb{W}(t),
\end{align}
where $\mathbb{C}=\mathbb{Q}^{-1}\mathbb{N}$. The noise in~\eqref{dY} has covariance matrix
\begin{align}\label{covar}
\mathbb{B}:=\mathbb{C}\mathbb{C}^{\top}=\mathbb{Q}^{-1}\mathbb{N}\mathbb{N}^{\top}(\mathbb{Q}^{-1})^{\top}.
\end{align}
Next we write
\begin{align}\label{newY}
\mathbb{Y}(t)=\mathbb{R}_{\omega t}\mathbb{Z}(t),
\end{align}
where $\mathbb{R}_s$ is the rotation~\eqref{rota}. Using the SDE~\eqref{dY} and It$\hat{o}$'s formula we obtain
\begin{align}\label{Ito}
d\mathbb{Z}(t)=-\lambda \mathbb{Z}(t)dt+\mathbb{R}_{\omega t} \mathbb{C} d \mathbb{W}(t).
\end{align}

Finally, in order to compare the process $\mathbb{Z}(t)$ with a standard two-dimensional Ornstein-Uhlenbeck process, we rescale time and space by writing
\begin{align}\label{rescale}
\mathbb{U}(t)=\frac{\sqrt{\lambda}}{\sigma} \mathbb{Z}(t/\lambda),
\end{align}
where
\begin{align}\label{sigma2}
\sigma^2=\frac{1}{2} {\rm tr}(\mathbb{B})=\frac{1}{2}(B_{11}+B_{22})=\frac{1}{2}\sum_{i,j=1}^2 C_{i,j}^2.
\end{align}

To identify the SDE that $\mathbb{U}(t)$ satisfies it is convenient to use the integrated form
\begin{align}\label{integrated}
\begin{split}
&\mathbb{U}(t)-\mathbb{U}(0)\\
&=\frac{\sqrt{\lambda}}{\sigma}\bigg(-\lambda \int_0^{t/\lambda} \mathbb{Z}(s)ds+\int_0^{t/\lambda} \mathbb{R}_{\omega s} \mathbb{C} d\mathbb{W}(s) \bigg)\\
&=\frac{\sqrt{\lambda}}{\sigma} \bigg(-\int_0^t \mathbb{Z}(u/\lambda)du+\frac{1}{\sqrt{\lambda}} \int_0^t \mathbb{R}_{\omega u/\lambda} \mathbb{C} d\tilde{\mathbb{W}}(u) \bigg),
\end{split}
\end{align}
where $\tilde{\mathbb{W}}(t)=\sqrt{\lambda}\;\mathbb{W}(t/\lambda)$ is a standard Brownian motion. Therefore
\begin{align}\label{dU}
d\mathbb{U}(t)=-\mathbb{U}(t)dt+\mathbb{R}_{\omega t/\lambda}\mathbb{D}d\tilde{\mathbb{W}}(t),
\end{align}
where $\mathbb{D}=(1/\sigma)\mathbb{C}$, so that ${\rm tr}(\mathbb{D}\mathbb{D}^{\top})=2$.

Theorem 1 of~\cite{Bax11} says that the distribution of $\{\mathbb{U}^{\lambda}(t): 0 \le t \le T \}$, where $\mathbb{U}^{\lambda}$ denotes the process satisfying~\eqref{dU} with a fixed $\lambda$, $\omega$, $\mathbb{U}^{\lambda}(0)=x$, fixed, and $0<T$ fixed, converges, as $\lambda /\omega$ goes to 0, to the distribution of the standard 2-dimensional O-U process $\{\mathbb{S}(t): 0 \le t \le T \}$ with $\mathbb{S}(0)=x$ and independent components. The proof uses the Martingale problem method (see~\cite{Bax11}).

The approximation~\eqref{Rotation} is then obtained by reversing the three changes of variables, starting from the SDE for the standard 2-dimensional O-U process
\begin{align}\label{OU2}
d\mathbb{S}(t)=-\mathbb{S}(t)dt+d\mathbb{W}(t).
\end{align}
In order to obtain an error bound on the approximation one would need a uniformity result with $T\rightarrow \infty$ as $\lambda /\omega \rightarrow 0$, in Theorem 1 of~\cite{Bax11}. In~\cite{GMW14} we applied this result to the system~\eqref{Kang3} and showed that the approximation~\eqref{Rotation} is quite good for the parameters we considered. The approximation predicts several critical quantities: the constant ratio of envelopes and phase difference of the E and I processes, their average instantaneous frequencies, and the relationship between the variance of instantaneous frequency and the amplitude for each of the two processes.  Simulations of the full model defined by~\eqref{Kang1} and~\eqref{Kang3}, using the same parameters used later in this paper, showed that these were all very close to the respective values predicted from the approximation. Impressively, the inverse relationship between instantaneous frequency variance and the envelope amplitude was just as predicted from the approximation~\eqref{db}, with large phase perturbations occurring more frequently when envelope amplitude was small. Moreover, as we show in the Results section, the approximation produces sample paths that qualitatively agree with those produced by the full model, as well as displaying the relationship between $d\phi$ and $Z(t)$ predicted by~\eqref{db}.

In terms of polar coordinates the approximation \eqref{Rotation} can be written as
\begin{equation}\label{polar}
\begin{split}
& \mathbb{V}(t) \approx \mathbb{V}^*(t)=\\
& \frac{\sigma}{\sqrt{\lambda}} \mathbb{Q} |\mathbb{S}(\lambda t)|  [\cos{(-\omega_d t+\phi(\lambda t))}, \sin{(-\omega_d t+\phi(\lambda t))}]^\top,
\end{split}
\end{equation}
where $ |\mathbb{S}(\lambda t)| =\sqrt{S_1(\lambda t)^2+S_2(\lambda t)^2}$ is the amplitude of $\mathbb{S}(\lambda t)$, and $\theta(t)=-\omega_d(t)+\phi(\lambda t)$ is the phase of $\mathbb{V}^*(t)$ at $\lambda t$, where $\phi(\lambda t) = \mbox{arg}[S_1(\lambda t)+iS_2(\lambda t)]$. 

For a standard two-dimensional Ornstein-Uhlenbeck process $\mathbb{S}(t)$, the stochastic process literature provides an explicit relationship between the process $Z(t)=|\mathbb{S}(t)|$, the \textit{radial process} or \textit{modulus process} associated with $\mathbb{S}(t)$, and $\phi(t)$, the \textit{phase perturbation process} (see \cite{Gard0x}, \cite{Borodin02}). The pair of processes, $Z(t)$ and $\phi(t)$, satisfies a system of stochastic differential equations
\begin{align}
dZ(t)&=\left[\frac{1}{2Z(t)}-Z(t)\right] dt + dW(t),\label{Z}\\
d\phi(t)&=\frac{1}{Z(t)} db(t),\label{db}
\end{align}
where $W(t)$ is a standard Brownian motion on the real line and $b(t)$ is a standard Brownian motion on the unit circle, independent of $W(t)$. We will see in the Results section that the approximation, Equation~\eqref{Rotation}, generates the same types of sample paths as the original model. Thus, from this point onward, we will take the approximation~\eqref{polar} as our model.

\section{Extension to a Population}

Suppose we have a family of $N$ E-I subpopulation models indexed by $i=1,2,...N$. We will have a version of \eqref{Kang1}, and hence of ~\eqref{polar}, for each member of the family, each with its own parameter values. The $\sigma_E, \sigma_I$  parameters for E and I subpopulations are all assumed to be equal for now, as are the time constants $\tau_E , \tau_I$. This can of course be changed to include distributions of parameters for the different subpopulations. These have been chosen to be biologically realistic insofar as possible. From~\eqref{polar}, \eqref{Z} and \eqref{db} we have, for each $i=1,2,...N$ (without coupling among subpopulations), for the phase of $\mathbb{V}^*_i(t)$,
\begin{align}\label{newdb}
d\theta_{i}(t)=-\omega_{di}dt+d\phi_i (\lambda_i t)
\end{align}
where 
\begin{align}\label{newslip}
d\phi_i(\lambda_i t)=\frac{1}{Z_i(\lambda_i  t)} db_i(\lambda_i t),
\end{align}
and for the amplitude of $V^*_i(t)$,
\begin{align}\label{newA}
\begin{split}
& dZ_i(\lambda_i t)=\\
& \frac{\sigma_i}{\sqrt{\lambda_i}}||\mathbb{Q}_i ||\left[\left(\frac{1}{2Z_i(\lambda_i t)}-Z_i(\lambda_i t)\right)dt +dW_i(\lambda_i t)\right],
\end{split}
\end{align}
where $||\mathbb{Q}||$ is the 2-norm of $\mathbb{Q}$.

An advantage of $\mathbb{V}^* (t)$ over $\mathbb{V}(t)$ is that when we transform to polar coordinates as in~\eqref{polar} we have the explicit expressions~\eqref{newdb} and~\eqref{newA} for the phase and amplitude processes, which are not available for $\mathbb{V}(t)$, and which are a natural way to describe quasi-cycle oscillators. This allows us to couple our family of quasi-cycle oscillators using Kuramoto-type coupling.

\section{Kuramoto Coupling}

The history of the Kuramoto model is beautifully set out in~\cite{Strogatz00}. See also~\cite{Acebron05} for a review. The model describes a dynamic system of phase changes that produces phase synchronization among a population of oscillators, with phase functions $\theta_j , j=1, 2, ...N$, depending on a distribution $g(\omega)$ of natural frequencies, $\omega_j$, and functions $\Gamma_{j,k}$ of phase differences:
\begin{align}\label{Kuramoto}
\frac{d\theta_j}{dt}=\omega_j+\sum_{k=1}^{N}\Gamma_{j,k}(\theta_k-\theta_j) \quad j=1,...,N.
\end{align}
The interaction functions $\Gamma_{j,k}$ could involve Fourier harmonics and an unspecified connection topology. In fact, Kuramoto studied the mean-field version of the model with sinusoidal coupling, 
\begin{align}\label{Kuramoto2}
\frac{d\theta_j}{dt}=\omega_j+\frac{K}{N}\sum_{k=1}^{N}\sin(\theta_k-\theta_j),\quad j=1,...,N.
\end{align}
For a collection of points $e^{i\theta_j}$ moving about on the unit circle in the complex plane, the quantity
\begin{align}\label{Krho}
\rho e^{i\psi}=\frac{1}{N}\sum_{j=1}^{N}e^{i\theta_j}
\end{align}
is the centroid of the phases, $\theta_j$, where the radius, $\rho(t)$, called the phase locking index, measures phase coherence, and $\psi(t)$ is the average phase. If the phases, $\theta_j, j=1...,N,$ move in a clump, the index $\rho(t)$ is near 1, whereas if the phases are spread out around the circle, the index $\rho(t)$ is near zero. One can write~\eqref{Krho} as
\begin{align*}
\rho e^{i(\psi-\theta_j)}=\frac{1}{N}\sum_{k=1}^{N}e^{i(\theta_k-\theta_j)}
\end{align*}
for any $\theta_j$. Equating the imaginary parts and substituting into~\eqref{Kuramoto2} one obtains
\begin{align*}
\frac{d\theta_j}{dt}=\omega_j+K\rho \sin(\psi-\theta_j),\quad~j=1, ..., N.
\end{align*}
In this form, we see that each phase $\theta_j$ is pulled toward the mean phase, $\psi$, with a strength proportional to the phase locking index $\rho$. If the population is becoming more coherent as time advances, $\rho(t)$ increases and so does $K \rho(t)$. Using the feedback loop, Kuramoto derived a critical value $K_c$ of $K$ such that for $K<K_c$, $\rho(t)\approx {\rm O}(1/\sqrt{N})$, whereas if $K>K_c$, $\rho(t)\approx \rho(K)+{\rm O}(1/\sqrt{N})$, where $\rho(K)$ is a saturating level of coherence~\cite{Strogatz00}.

\subsection{Kuramoto Coupling of a Family of Quasi-cycle Oscillators }

Results for stochastic versions of the Kuramoto model have been obtained, e.g., by~\cite{Sonn13, deville12}. Our question here is whether these results for the Kuramoto model~\eqref{Kuramoto2}, in particular monotonicity of coherence with coupling strength and rapid change of coherence near a critical value of coupling strength, will continue to hold for a coupled version of the stochastic systems~\eqref{newdb} and \eqref{newA} for $i=1, ..., N$. 

We modified the system \eqref{newdb} and \eqref{newA} by inserting a coupling term like that in~\eqref{Kuramoto} into their drift coefficients, assuming Theorem 1 of~\cite{Bax11} holds uniformly over $i=1,...,N$. The factor $Z_j/Z_i$, the ratio of amplitudes, appears in the coupling of phases as in~\cite{Daff11, Ton14}. It reconciles the fact that subpopulations $i$ and $j$ have different amplitude processes, and their ratio is a separate factor from the coupling strength determining the impact of the phase subpopulation $j$ on the phase of subpopulation $i$. This point is pursued in the Discussion. The result is a system of coupled phases, 
\begin{align}\label{GMW1}
\begin{split}
& d\theta_{i}(t)=\\
& \left(-\omega_{di}+\frac{1}{2N}\sum_{j=1}^{N} \frac{Z_j(\lambda_jt)}{Z_i(\lambda_it)} \mathbb{C}_{ij} \sin(\theta_{j}(t)-\theta_{i}(t))\right)dt\\
& +\frac{db(\lambda_it)}{Z_i(\lambda_i t)},
\end{split}
\end{align}
where $\mathbb{C}$ is a matrix defining the coupling between subpopulations, $\mathbb{C}_{ij}\geq0, i\neq j; \mathbb{C}_{ii}=0,  i=1,...N$. In addition to impacting each others' phases, we would expect that the amplitudes of the various E-I subpopulations would impact each other as well. This is important because the amplitudes affect the phase slip processes~\eqref{newslip}. As in an OU process, where the process is pulled toward 0 in proportion to its distance from 0, we can conceptualize the driving impact of each subpopulation's amplitude on that of the others as  a `pull' on the others' amplitudes proportional to the difference between them. This symmetrical coupling has the advantage that it doesn't allow the amplitudes to become too large, which is necessary in order that the phase slip processes~\eqref{newslip} be observable. We recognize that other forms of amplitude coupling are possible and that some of these might be more realistic in neural systems. We come back to this point in the Discussion. We couple the amplitudes \eqref{newA} as follows:
\begin{align}\label{GMW2}
\begin{split}
& dZ_{i}(\lambda_i t)=\\
& \left[\frac{\sigma_i}{\sqrt{\lambda_i}}||\mathbb{Q}_i ||\left(\frac{1}{2Z_{i}(\lambda_i t)}-Z_{i}(\lambda_i t)\right)\right] dt\\
& +\left[\frac{1}{2N} \sum_{j=1}^{N} \mathbb{C}_{ij} (Z_{j}(\lambda_j t)-Z_i(\lambda_it))\right] dt\\
& +\frac{\sigma_i}{\sqrt{\lambda_i}}||\mathbb{Q}_i ||dW_{t}(\lambda_i t).
\end{split}
\end{align}

The overall degree of coupling in the entire set of subpopulations can be characterized in general by the 2-norm of the matrix $\mathbb{C}$, $||\mathbb{C}||$, which when normalized with regard to the population size, $N$, can be seen to play the same role as the coupling coefficient, $K$, in the Kuramoto model. 

\section{Results}

In order to study the consequences of this Kuramoto-type coupling of quasi-cycle oscillators, we solved \eqref{GMW1} and \eqref{GMW2} numerically using the Euler-Maruyama method \cite{Kloeden} for a range of the number of coupled E-I systems,  $N$, and  the overall coupling strength between the systems, $||\mathbb{C}||$, similar to Kuramoto's $K$. We assumed that the subpopulations all oscillated in a narrow band of natural frequencies distributed normally with mean 437.72 rad/s and standard deviation 1 rad/s, similar to the standard system studied for the Kuramoto model. We also began each realization with the phases of the various systems distributed uniformly between $-\pi$ and $\pi$. Finally we studied all-to-all coupling matrices $\mathbb{C}_{i,j}$, all equal except that $C_{i,i}=0, \forall i$. We investigated the coherence measure, $\rho$, as a function of $||\mathbb{C}||$, and also studied the dynamics of the coherence by looking at the membership of the synchronous and asynchronous E-I groups via raster plots and plots of the mean frequency, amplitude, and phase in the most synchronous group compared with those of the entire set of processes for representative values of $||\mathbb{C}||$.

In our simulations we chose a natural, or resonant, frequency,  $\omega_d$, that for different E-I processes was chosen from a (clipped, i.e. compact support)  Gaussian distribution. The damping rate, $\lambda$ and the noise amplitude, $\sigma$, also were required to vary in a consistent way because both quantities depend, in different ways, upon the same parameters from \eqref{Kang1}, namely the synaptic efficacies $S_{i,j}, i,j=E,I$, the time constants $\tau_E, \tau_I$, and $\sigma_E, \sigma_I$ (see \eqref{Damping_rate},~\eqref{omega_d},~\eqref{sigma}). For purposes of the present study, $\lambda$ was computed for each separate process from the following equations given $\omega_d$ for that process (many significant digits retained because of the sensitivity of the calculations),
 \begin{align*}
 \begin{split}
 & S_{II}=0.5-\\
 & \frac{\sqrt{(-27778.33)^2 - 4(194443.33-\omega_d^2) (6944.44)}}{2(6944.44)},
 \end{split}
\end{align*}
 \begin{align*}
\lambda=83.33 S_{II},
\end{align*}
under the assumption that $S_{EE}=1.5, S_{IE}=4, S_{EI}=1, \tau_E=0.003, \tau_I=0.006$, $\sigma_E=\sigma_I=12$. Thus, the variation in $\omega_d$ was attributed entirely to variability in $S_{II}$ (see~\cite{GMW14} for a discussion of the roles of synaptic efficacies in this model). Similarly, $\sigma$ turns out to be
\begin{equation*}
\sigma=2998.38/\omega_d,
\end{equation*}
and $\mathbb{Q}$ is
\begin{equation*}
\mathbb{Q}=\begin{pmatrix}
			-\omega_d & \lambda+166.67\\
			0 & 666.67
			\end{pmatrix}.
\end{equation*}
For the parameter values we used in our simulations of~\eqref{GMW2}, $||\mathbb{Q}|| \approx 703.5$. We rescaled the amplitudes into a range that allows for an observable phase slip process by using $||\mathbb{Q}||/703.5\approx1$ instead of $||\mathbb{Q}||$ in those simulations. As mentioned earlier, if the amplitudes become too large, the phase slip processes ~\eqref{newslip} become negligible and we have what amounts to a standard Kuramoto model.

\begin{table}[h!]
\caption{Parameters used in simulations and for figures.}\label{Table1}
\begin{center}
\begin{tabular}{|c|c|c|}
\hline
Variable & Mean Value & Units\\
\hline
\hline
$\lambda$ & 8.333 & 1/seconds\\
\hline 
$\omega_d$ & 437.72  or 69.66 & radians per second or Hz\\
\hline
$\lambda/\omega_d$ &$0.019$ & dimensionless\\
\hline
$\sigma$ & $6.85$ & mV\\
\hline
$\Delta t$ & 0.00005 & seconds\\
\hline
\end{tabular}
\end{center}
\end{table}

\begin{figure}[!ht]
\begin{center}
\includegraphics[width=3.3in]{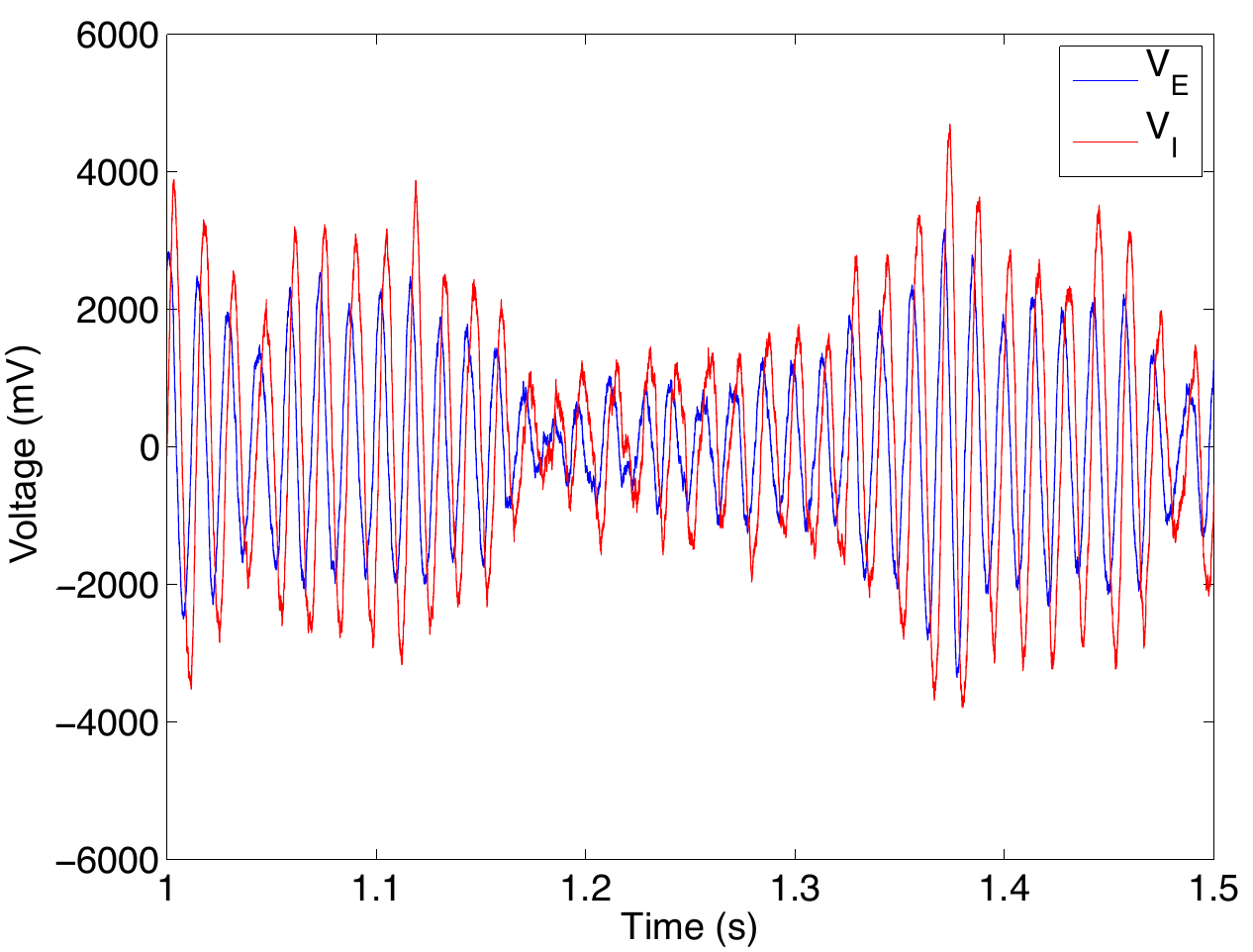} 
\end{center}
\caption{Sample path from Equation \eqref{Rotation}.The natural frequency, $\omega$, of the E-I process was set to 69.66 Hz (437.72 rad/s). Other parameters as in Table 1.} 
\label{Figure1}
\end{figure}

\subsection{Uncoupled oscillations}

\begin{figure}[!ht]
\begin{center}
\includegraphics[width=3.3in]{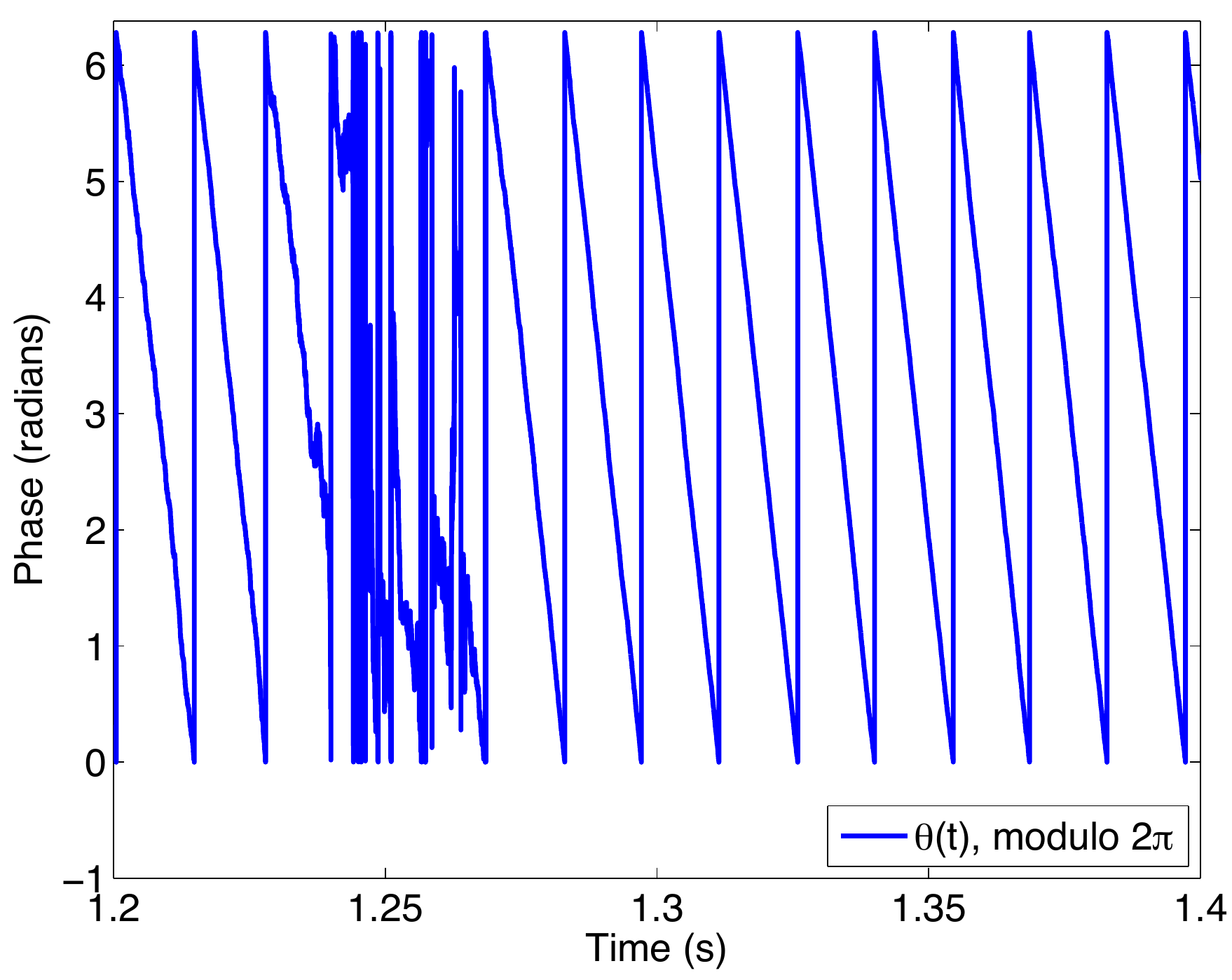}  
\end{center}
\caption{Single phase process from \eqref{GMW1}, in the  the absence of coupling (i.e., $\mathbb{C}=\mathbb{O}$, the zero matrix). The data are from a realization using the same parameters as in Figure \ref{Figure1}.}
\label{Figure2}
\end{figure}

\begin{figure}[!ht]
\begin{center}
\includegraphics[width=3.3in]{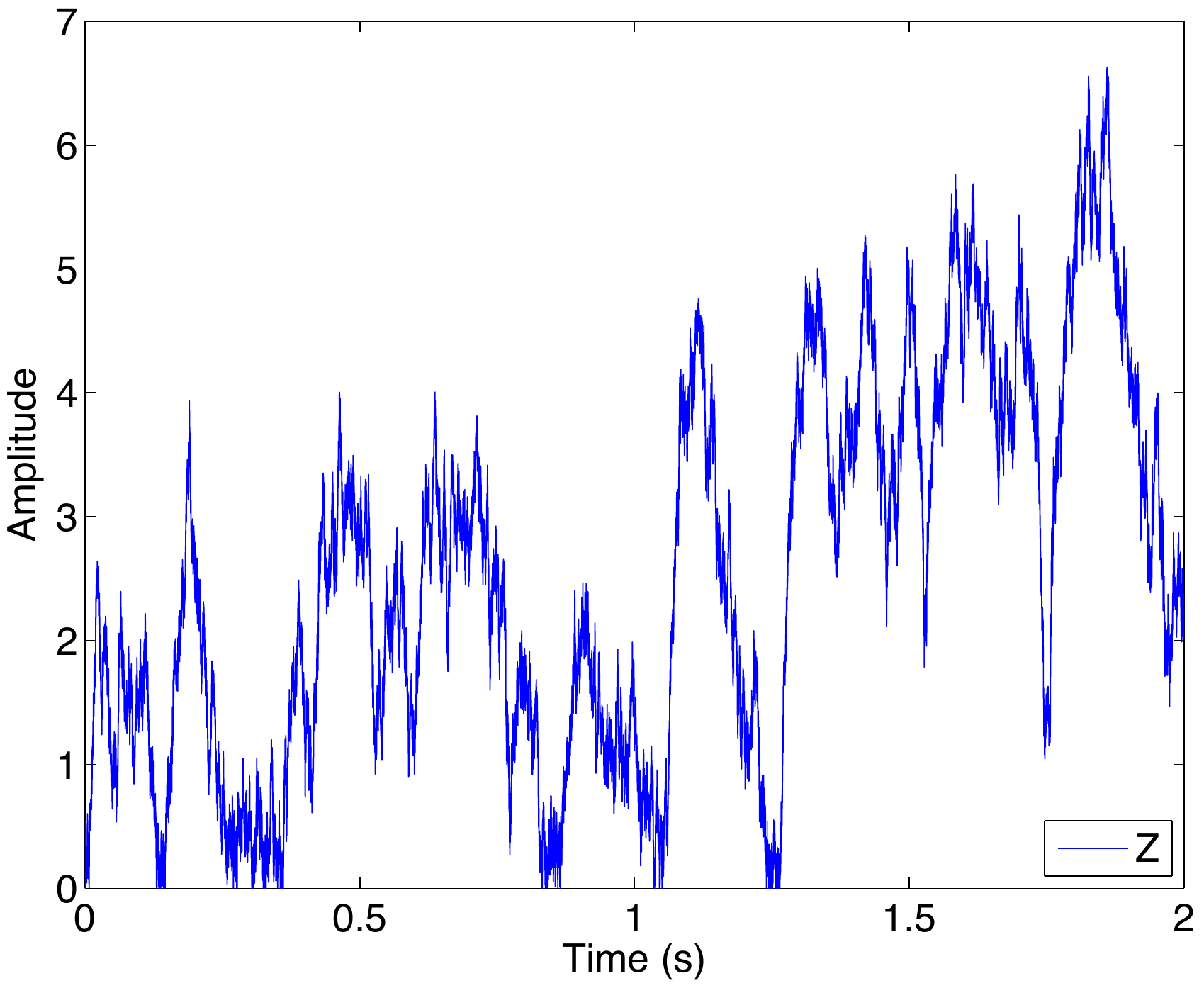} 
\end{center}
\caption{Single amplitude process, $Z(t)$, from \eqref{GMW2}, in the  the absence of coupling (i.e., $\mathbb{C}=\mathbb{O}$, the zero matrix). The data are from the same realization as in Figure \ref{Figure2}.} 
\label{Figure3}
\end{figure}

Before we describe our results for the coupled systems defined by \eqref{GMW1} and \eqref{GMW2}, we need to establish that the phase and amplitude processes, \eqref{newdb} and \eqref{newA} produce appropriate stochastic oscillatory paths. To this end we first simulated \eqref{Rotation} to demonstrate the ability of the approximation $\mathbb{V}^*(t)$ to generate oscillatory sample paths qualitatively similar to those generated by the full model \eqref{Kang1}. One such sample path is shown in Figure \ref{Figure1} and it is qualitatively the same as the sample paths produced by the full model, including the gamma bursts (episodes of high amplitude oscillations, cf. Figure 6 in~\cite{GMW14}) produced by that model. We then simulated \eqref{GMW1} and \eqref{GMW2} in the absence of coupling ($\mathbb{C}=\mathbb{O}$ (the zero matrix)), which are in fact the same system as \eqref{newdb} and \eqref{newA}, and inspected the paths of the individual uncoupled phase and amplitude processes. Figure \ref{Figure2} displays a typical phase path, and Figure \ref{Figure3} a typical amplitude path for a single E-I process approximated by \eqref{GMW1} and \eqref{GMW2}. The phase and amplitude processes of Figures \ref{Figure2} and \ref{Figure3} are not those of the model (6) of \cite{GMW14}, which is \eqref{Kang1} of the present paper, but are samples from the almost equivalent model $V^* (t)$ of~\eqref{polar}. The oscillatory nature (rotation in polar coordinates) of this process is revealed by the orderly phase progression. The amplitude process can be regarded as the envelope of the oscillations.

\begin{figure}[!ht]
\begin{center}
\includegraphics[width=3.3in]{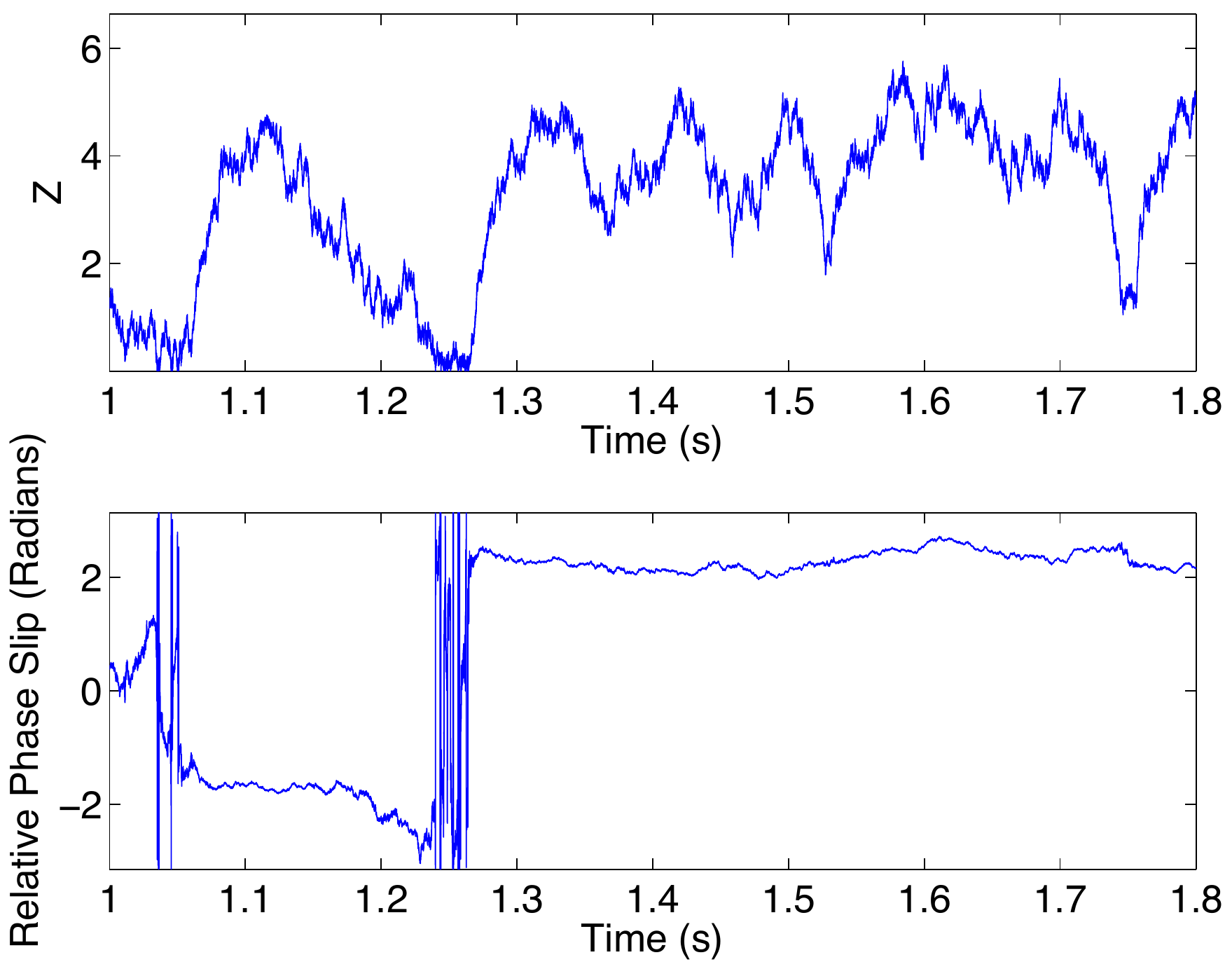}
\end{center}
\caption{Amplitude, $Z(t)$ and phase, $\phi(t)$, of a single process plotted together. The data are from the same realization as that in Figure \ref{Figure2} and Figure \ref{Figure3}. Notice that large phase slips occur when the amplitude is low, around 1.03 s and 1.25 s in this figure.}
\label{Figure4}
\end{figure}

Note in Figure \ref{Figure2}, however, that the orderly phase progression is interrupted on occasion by large increments of the phase perturbation process. These increments are related to the amplitude process as predicted by Equation \eqref{db}. An example of this is shown in Figure \ref{Figure4} where phase is plotted with amplitude in the same realization. Large phase perturbation increments occur when amplitude is low and small increments occur when amplitude is large.

\subsection{Phase coherence and coupling}

For the purposes of this paper we define synchronization in terms of phase coherence, which is the Kuramoto Phase Locking Index (PLI) in \eqref{Krho}, 
\begin{align*}
\rho=\frac{1}{N}\left|\sum_{j=1}^{N}e^{i\theta_{j}}\right|.
\end{align*}
The $\theta_{j}$ are the phases of the $N$ separate stochastic processes $V_j^*(t)$.

\begin{figure}[!ht]
\begin{center}
\includegraphics[width=3.3in]{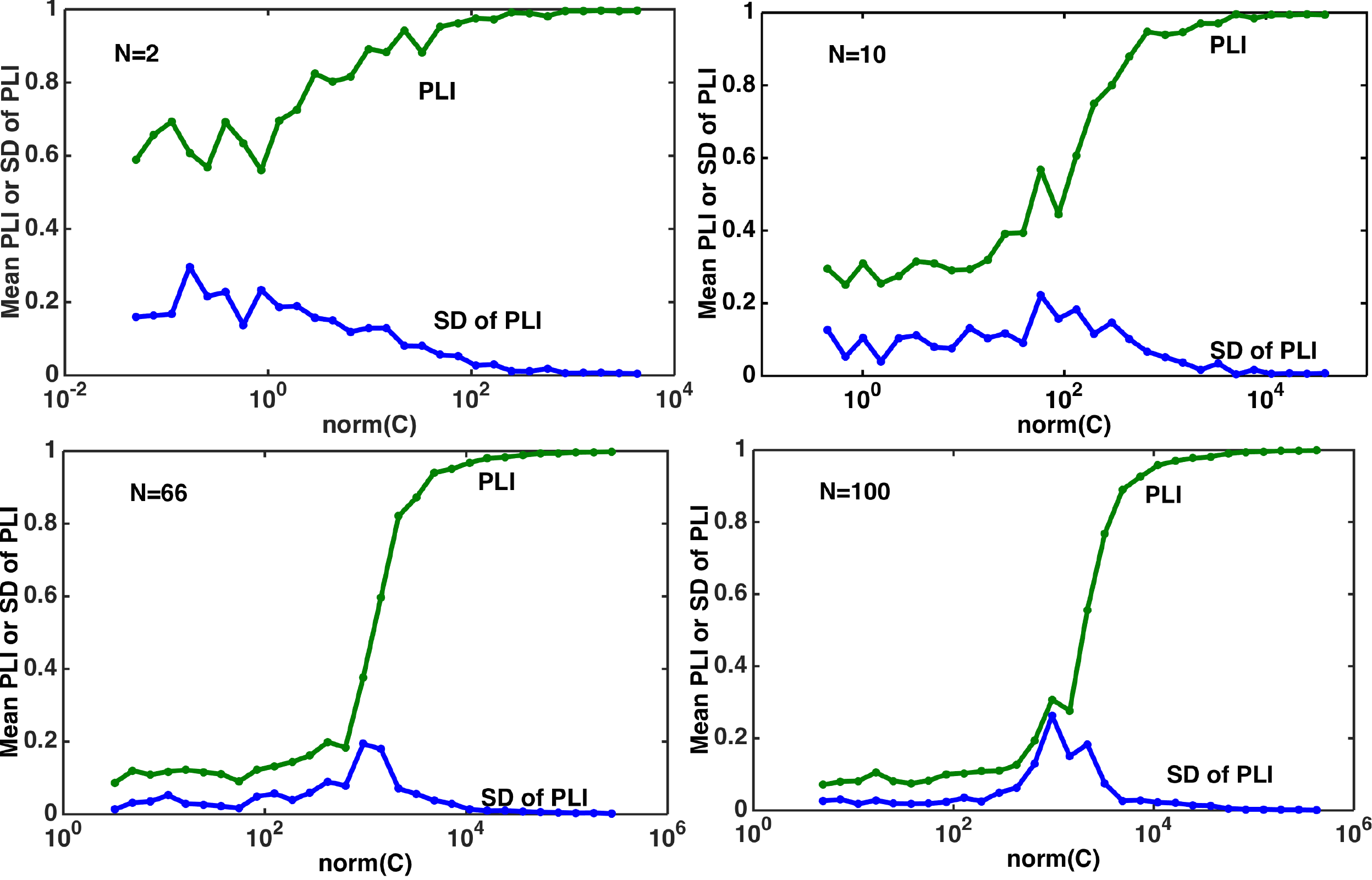}
\end{center}
\caption{Phase coherence, $\rho$ (phase-locking index, PLI), as a function of coupling strength, $||\mathbb{C}||$, with normally distributed natural frequencies. In all graphs, the connected points represent means over 10 independent realizations of $\bar{\rho}$, and the solid lines at the bottom of the graphs are the standard deviations of $\bar{\rho}$ over the 10 realizations. For each realization, $\bar{\rho}$ entered into the interation average was the average value of $\rho$ over 5000 iterations (0.25 s) of the model, beginning after iteration 5000, at which point the model was at or close to a steady state. Natural frequencies, $\omega_d$,  of E-I processes were sampled from a (clipped) Gaussian distribution with a mean of 437.72 rad/s and a standard deviation of 1 rad/s, and this was done anew for each iteration and each coupling value. Starting phases were randomly distributed between $-\pi$ and $\pi$, and starting amplitudes were uniformly distributed between 0 and 1.}
\label{Figure5}
\end{figure}

Figure \ref{Figure5} displays the average phase coherence of 10 independent realizations of \eqref{GMW1} and \eqref{GMW2} as a function of coupling strength, $||\mathbb{C}||$, among E-I populations with $N$=2, 10, 66, and 100,  with (clipped between -3 and +3) normally-distributed natural frequencies. As in the basic and stochastic Kuramoto models, phase coherence is nearly constant for coupling strengths below a critical value, and increases rapidly (for $N\geq 10$) with increasing coupling strength until near an asymptote at which the phases are all locked together and $\rho$ stays very near 1, with only occasional slight departures caused by occasional large noise samples in one or more of the coupled processes. The standard deviation of $\rho$ over the 10 realizations in each graph is roughly constant for small coupling values except near the critical coupling value, where phase coherence is rising rapidly under the mutual influence of the different E-I processes. The standard deviation of $\rho$, of course, decreases dramatically after PLI reaches its asymptote. The effect of increasing the number of coupled systems is to increase the critical coupling value, in terms of $||\mathbb{C}||$, and to decrease the minimum phase coherence of the system. Clearly stronger coupling is required to synchronize a larger number of independent systems. 

\begin{figure}[!ht]
\begin{center}
\includegraphics[width=3.3in]{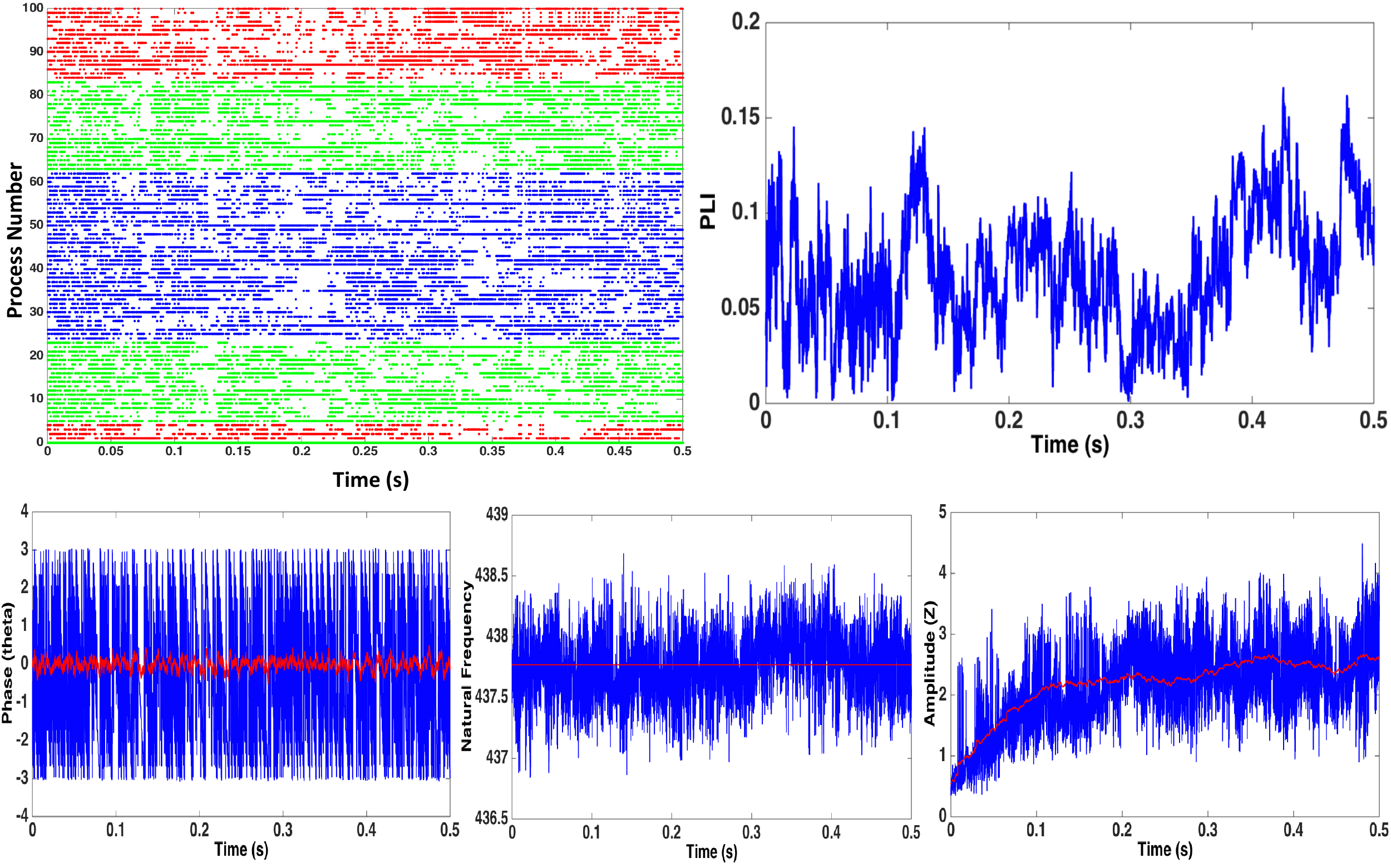}
\end{center}
\caption{Uncoupled processes - minimal synchronization. No coupling between processes, single
realization with $N = 100$, $||\mathbb{C}_{i,j}|| = 0$. In the raster plot in the upper left each plotted point
represents membership of the indicated process in the `synchronous group' (see text) at the
indicated iteration. Red points represent processes whose natural frequency fell within the extreme
two highest and two lowest of 20 equal-sized bins over the range of the normal distribution of
natural frequencies; blue points represent processes whose natural frequencies fell within the two
most populous bins; green points represent the remainder of processes, with natural frequencies
between the extremes and the most populous bins. Upper right plot: PLI over the 100 processes
at each time point. Lower left plot: Red line indicates mean phase (theta) of all 100 processes;
blue line indicates mean phase of processes in the synchronous group. Lower middle plot: Red line
indicates mean natural frequency of all 100 processes; blue line indicates mean natural frequency
of processes in synchronous group. Note that instantaneous frequencies are not always equal to
natural frequencies because of phase slips caused by driving noise. Lower right plot: Red line
indicates mean amplitude of all 100 processes; blue line indicates mean amplitude of the processes
in the synchronous group.}
\label{Figure6}
\end{figure}

\subsection{Dynamics of phase coherence}
\begin{figure}[!ht]
\begin{center}
\includegraphics[width=3.3in]{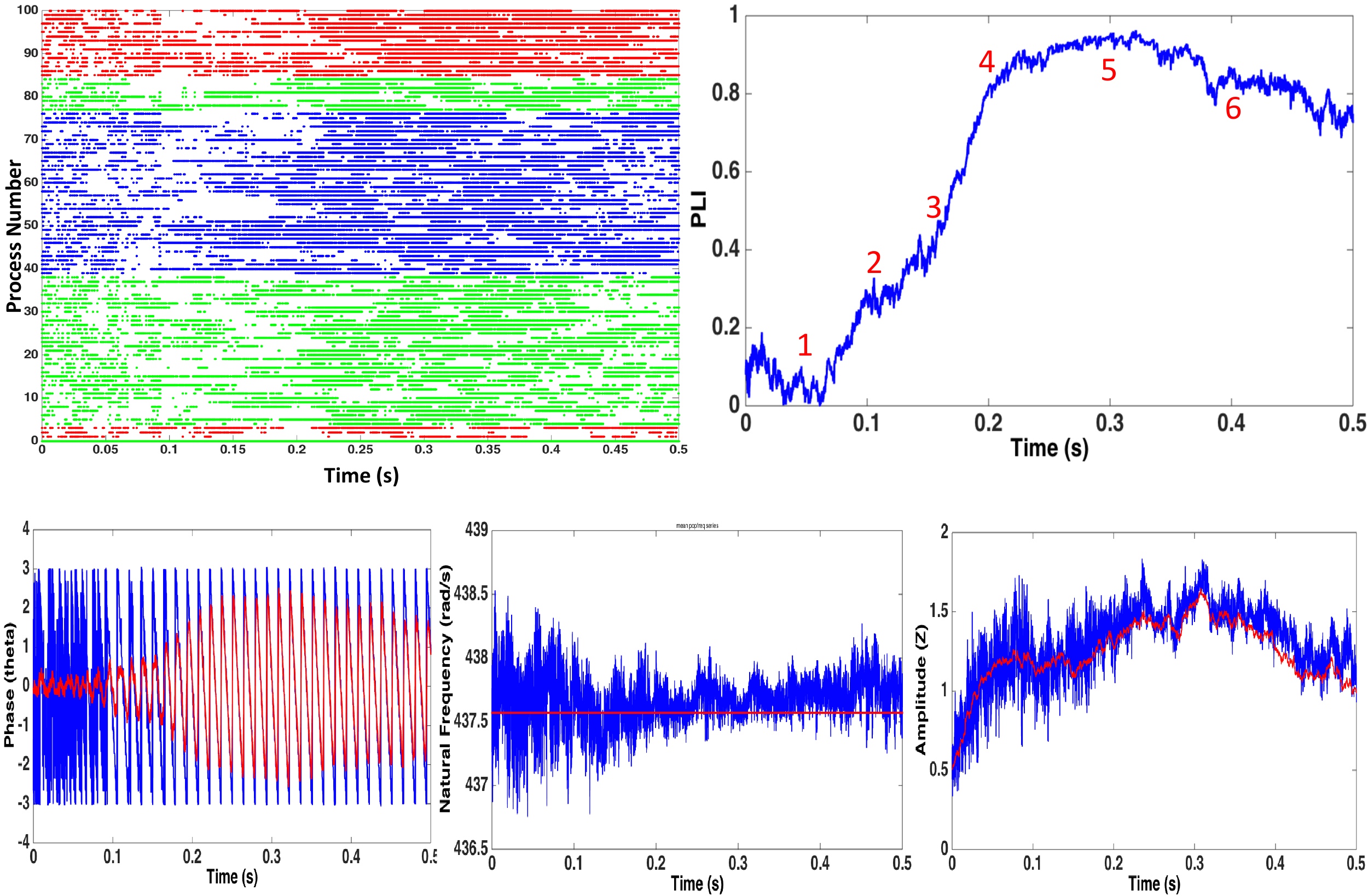}
\end{center}
\caption{Coupled processes - substantial synchronization. Relatively strong coupling between processes, single
realization with $N = 100$, $||\mathbb{C}_{i,j}|| = 4950$. In the raster plot in the upper left each plotted point
represents membership of the indicated process in the `synchronous group' (see text) at the
indicated iteration. Red points represent processes whose natural frequency fell within the extreme
two highest and two lowest of 20 equal-sized bins over the range of the normal distribution of
natural frequencies; blue points represent processes whose natural frequencies fell within the two
most populous bins; green points represent the remainder of processes, with natural frequencies
between the extremes and the most populous bins. Upper right plot: PLI over the 100 processes
at each time point. Lower left plot: Red line indicates mean phase (theta) of all 100 processes;
blue line indicates mean phase of processes in the synchronous group. Lower middle plot: Red line
indicates mean natural frequency of all 100 processes; blue line indicates mean natural frequency
of processes in synchronous group. Note that instantaneous frequencies are not always equal to
natural frequencies because of phase slips caused by driving noise. Lower right plot: Red line
indicates mean amplitude of all 100 processes; blue line indicates mean amplitude of the processes
in the synchronous group.}
\label{Figure7}
\end{figure}

How do the group dynamics of the unsynchronized and synchronized collections of the stochastic E-I pairs evolve, and how do they compare? To gain insight, we determined, using simulations similar to those in Figure \ref{Figure5}, for two coupling values, $||\mathbb{C}||=0$ and $||\mathbb{C}||=4950$, and with $N=100$, which particular pairs of E-I processes, characterized by their natural frequencies, participated in the group of processes that was most highly synchronized. To do this we sorted the phases of the processes into 20 equal phase bins between $-\pi$ and $\pi$ and chose the one with the largest membership at each iteration as the `synchronous group' at that moment. Figures \ref{Figure6} and \ref{Figure7} display examples of raster plots indicating the membership of the various processes in this synchronous group for the two different coupling strengths. The processes are displayed in the plot in ascending order of their natural frequency, $\omega_j$. Figures \ref{Figure6} and \ref{Figure7} also display plots of phase locking index, PLI, $\rho$, computed across all 100 processes as a function of time since the beginning of the simulation. Included are plots of the mean phase, natural frequency, and amplitude of the members of the synchronous group at each time point of the simulation compared with those of the total set of 100 processes.

The raster plots of synchronous group membership clearly show that when coupling is 0 (Fig. \ref{Figure6}), membership is highly variable and the synchronous group itself is sporadically populated, with most of the processes involved at some point regardless of natural frequency. There is no strong tendency for either extreme or central frequencies to dominate. For stronger coupling (Fig. \ref{Figure7}), the result is similar, except that more processes are in the synchronous group more of the time, and the processes with higher natural frequencies dominate the synchronous group. Apparently the independent noises for phase and amplitude for each process preclude most processes from settling into the synchronous group for a very long time, although a few do seem to be in it most of the time, including especially those with higher natural frequencies which would tend to lead the others in phase. 

In light of the heterogeneity of synchronous group membership across model time, the plots of mean frequency and mean amplitude of the processes in the synchronous group with stronger coupling are interesting. Figure \ref{Figure7} reveals that, on average, when coupling is relatively strong the mean natural frequency of the synchronous group is higher than that of the entire population. Figure~\ref{Figure7} also shows that the mean amplitude of the processes in the synchronous group is higher than the mean amplitude of the entire population. It is reasonable that processes with higher amplitudes and frequencies would tend to dominate the interactions.

\begin{figure}[!ht]
\begin{center}
\includegraphics[width=3.45in]{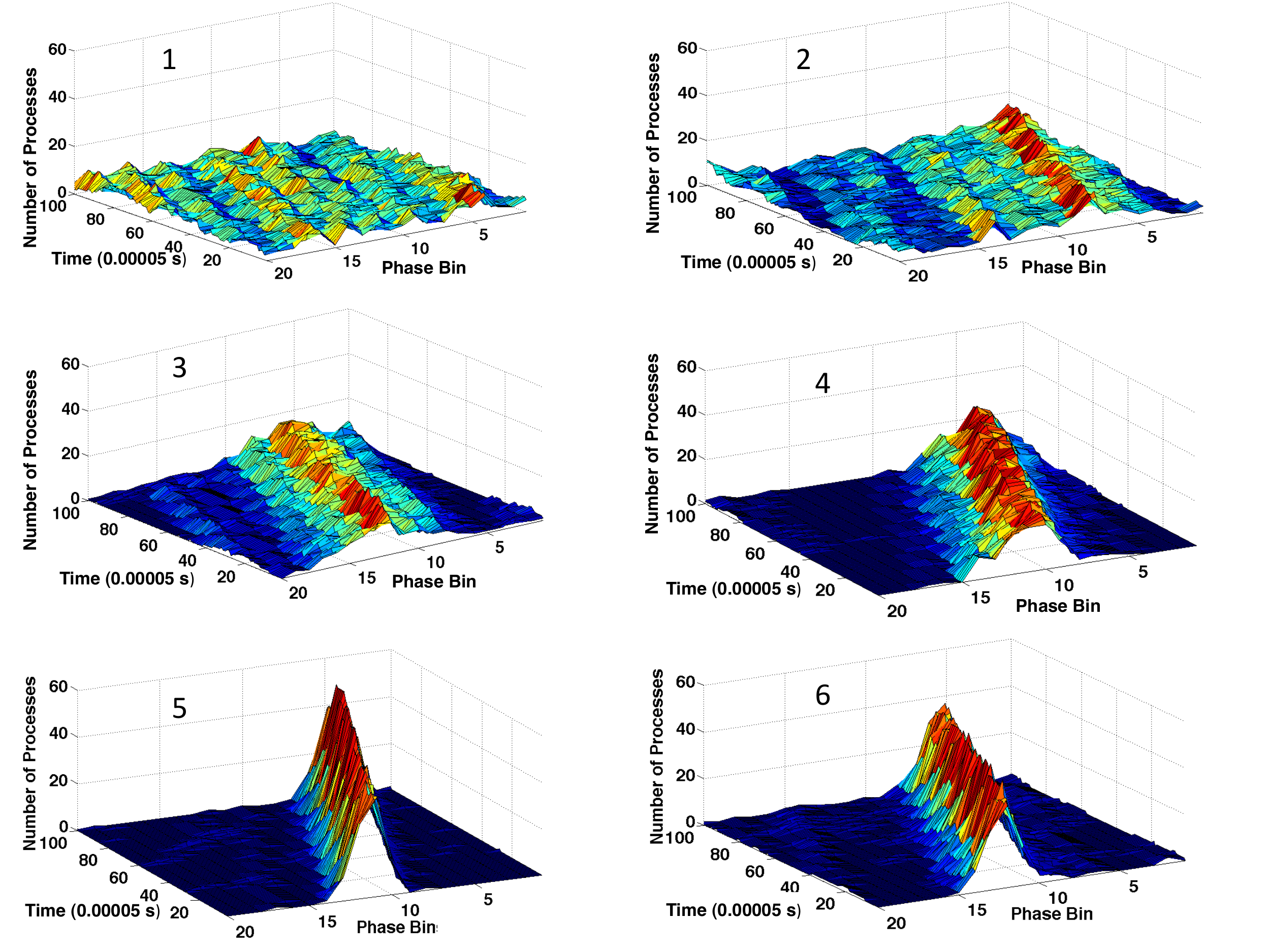}
\end{center}
\caption{Dynamics of synchronization. Distribution of phases across 20 equal-sized phase bins
from $-\pi$ to $\pi$ at different time points in a single realization (same as in Figure \ref{Figure7}) of the model for
relatively strong coupling between processes, showing the progression of the neural oscillators
from an unsynchronized to a synchronized state. Each plot represents a 5 ms period of the
realization. Plots from upper left to lower right progress from early to later in the realization;
specific time points are indicated by numbers that match the red numbers in the upper right plot
in Figure \ref{Figure7}, which also plots the associated coherence values.}
\label{Figure8}
\end{figure}

Figure~\ref{Figure8} displays additional detail concerning the evolution of synchronization in the case of the realization shown in Figure~\ref{Figure7}. The number of members in each phase bin is plotted in Figure~\ref{Figure8} for a series of 0.00005 second time slices. The `synchronous group' is represented by the highest peak of the distribution at each time point. Clearly, early in the realization, the phases of the various processes are incoherent, spread across the 20 bins roughly equally as they were at the onset of the realization. Gradually, as model time progresses, the interactions between the processes promoted by the moderately strong coupling creates first several clumps of synchronized processes, and then a single dominant group. Notice that, consistent with the slight drop in PLI later in the realization, at around 0.4 seconds (labeled `6' in the graph of PLI in Fig.~\ref{Figure7}), the dominant group becomes somewhat less coherent again, presumably under the influence of the noise that is driving the quasi-cycles. 

\section{Discussion}

We have reported results from a combination of a linear stochastic model of populations of neurons with the Kuramoto approach to synchronization of quasi-cycle oscillators. This is new because for the first time we demonstrate for quasi-cycles the existence, and study the properties, of a critical coupling leading to their synchronization. Such oscillators are likely to be representative of a broad class of biological oscillators, such as neurons and neural networks, and including gene networks. The E-I processes in our model could represent a range of possible neurological configurations, from pairs of reciprocally-interacting neurons, through populations of sparsely-firing neurons that generate an oscillating field potential, to synchronously firing populations of neurons.

We have demonstrated that for this novel quasi-cycle stochastic Kuramoto model for finite $N$, as for other Kuramoto models, phase coherence generally increases with coupling strength. The incoherent state at low coupling changes first to a partially coherent state and then to a nearly phase-locked state as coupling increases beyond a critical value. We have also shown some of the details of the transition in the form of raster plots and other properties of the group of most nearly phase-locked oscillators. The exploration of additional coupling matrices is left to future work. If quasi-cycles are indeed the source of neural oscillations, as argued in~\cite{GMW14}, then the results in the present paper become relevant to synchronization of neural systems in real brains.

In the present work we introduced, in addition to phase coupling, coupling of the amplitudes of the quasi-cycle oscillators. This was necessary because the amplitudes of the oscillators affect their phase slip processes, and thus also affect the synchronization of the phases. If phase slips are large and frequent (i.e., when amplitudes are very small), synchronization becomes more problematic, whereas if phase slips are small and infrequent (i.e., when amplitudes are very large) the rotation of each quasi-cycle is more regular and the situation more closely resembles interactions among a set of pure phase oscillators. One limitation of our present work is that the form of the amplitude coupling we used was symmetrical, i.e., the lower amplitude oscillators pulled the higher ones down, and vice versa, so that the amplitudes tended to stabilize in a small range that allowed for observable phase slips. We also scaled the amplitudes to maintain this range. It is possible that other forms of amplitude coupling would operate similarly. For example, we tested a model with the ratio of amplitudes, $Z_j(\lambda_j t)/Z_i(\lambda_i t)$ inserted in place of the difference of amplitudes in~\eqref{GMW2}. Results were highly similar, although the amplitudes tended to stabilize at somewhat higher levels. On the other hand, inserting only the amplitude of the `other' coupled oscillator, $Z_j(\lambda_j t)$, in place of its difference from $Z_i(\lambda_i t)$ results in an explosion of amplitudes and thus an exit from the quasi-cycle regime. Similarly, a coupling factor of $cos(\theta_j-\theta_i)$ or $Z_j(\lambda_j t) cos(\theta_j-\theta_i)$ also has this effect. Thus, these couplings will not work for quasi-cycle oscillators.

All of this raises the question of what a biologically-realistic coupling of quasi-cycle oscillators would be, for example in the context of a neural network. Typically, neurons are coupled via synapses. Our original model~\eqref{Kang1} expresses the coupling between E and I neurons realistically as a bidirectional, short-range synaptic coupling. But when we move into the realm of the quasi-cycle oscillator we lose the ability to represent directly the long-range couplings between the neurons themselves, and can only represent couplings between the quasi-cycle oscillators' phases and amplitudes. It remains to be worked out how the couplings we have studied are related to realistic couplings between neurons or other biological entities. Long range couplings between neurons are more common between excitatory neurons, but can also occur between I and E neurons. Thus, the coupling we studied in the present work might best capture a network of quasi-cycle oscillators in which bidirectional E and I connections occur at two scales, say within micro-columns (the individual E-I quasi-cycle oscillators) and also between columns (the coupled set of quasi-cycle oscillators). Such a network would allow for both amplitude and phase `pulls' in both directions, as in~\eqref{GMW1} and \eqref{GMW2}. Our conclusion would then be that the larger-scale set of oscillators synchronizes according to coupling strength in a Kuramoto-like manner. 

\subsection{Relation to other models}

Although we have not provided mathematical results for the critical coupling of quasi-cycles, there do exist such results for related models. For example, the paper of Ott and Antonsen~\cite{Ott08} considers what happens when a system of coupled oscillators is also coupled to an external oscillator. It is found, using an Ansatz regarding the form of the coefficients in the Fourier expansion of the Fokker-Planck (Kolmogorov) equation for the stochastic Kuramoto system, that in the large $N$ limit there is a critical value of the coupling constant at which there is a phase transition from incoherence to partial or complete synchronization. 

Sonnenschein and Shimansky-Geier~\cite{Sonn13} derived a similar result for a stochastic model with Gaussian noise, assuming that the phases of the oscillators are Gaussian in distribution. In their model, the variance of the phases of the oscillators, $\sigma ^2 (t)$ in their notation, goes to infinity as $t\rightarrow \infty$ for $C<2D$, and has a finite limit for $C>2D$, where $C$ is the coupling strength and $D$ is the variance of the noise. It is quite possible that a similar balance between coupling strength and noise variance holds in our model, \eqref{GMW1} and \eqref{GMW2}, although our noise variance in \eqref{GMW1} involves $Z_i$. Moreover, $2D$ in their model is a transition point between two regimes of phase variance, whereas in our model there is a critical coupling strength that separates incoherence from coherence. Nonetheless, if we consider a normalized coupling matrix 2-norm, $|| \mathbb{C}||/(N^2 -N)$, as the coupling strength, and the $N=100$ case in which $\rho$ is near 1 at $|| \mathbb{C}|| \approx 10^4$, then the critical coupling value is around $10^4/(10^4 -100) \approx 1$ in our model. Because $Z_i(\lambda_i t)$ converges to values around 1 when coupling is large (or even moderate, see Figure \ref{Figure7}), and the variance of the noise is 1 in our simulations, at convergence an analogue of the $D$ of~\cite{Sonn13}, $(db(\lambda_i t)/Z_i(\lambda_i t))^2$, is approximately equal to $1/1 =1$. Their theory would predict a value of $2 \times D = 2 \times 1 = 2$,  for our critical coupling value, which is of the same order of magnitude as what we observed, although this could certainly change as $N$ becomes even larger. 

It should be noted, however, that our approximation informs us  of an additional constraint on the relation between critical coupling value and the noise that is driving the quasi-cycles. Namely, $\sigma/\sqrt{\lambda}$ must not go to either 0 or infinity. Because $\lambda$ must be small relative to $\omega$, $\sigma$ too must be relatively small in order for the quasi-cycles to be of at least moderate size. Thus, because the noise is constrained, if there is a relationship between the critical coupling value and the noise, that relationship is also constrained. 

The recent paper of Ton et al.~\cite{Ton14} seems to be one of the closest to ours in some
respects. The deterministic dynamics, however, is taken from a Freeman neural
mass model expressed by a second order ODE rather than from a linear stochastic model that generates quasi-cycles; E-I pairs are involved but in a different way. The emphasis is on the effect
of delays, which are transformed into phase shifts. Deterministic limit cycles around
an unstable focus are the source of unstable oscillations, rather than noisy quasi-cycles near a stable focus as in the present paper. A commonality is that the limit
cycles in~\cite{Ton14} have a variety of radii and each coupling parameter contains a factor of
the ratio of these radii, just as in~\cite{Daff11}, with justification involving, again, averaging
over one period of the oscillation (see~\cite{Ton14}, Methods). We follow this procedure in
introducing the ratio of amplitudes into~\eqref{GMW1}, where the distances of the quasi-cycle oscillations from the fixed point vary from moment to monent for each oscillator as well as across oscillators at any given moment. The introduction of such a factor into the phase coupling is not critical in our model, however, for we found similar results when it was omitted and the coupling was exactly as in the Kuramoto model.

The object of study in~\cite{Ton14} is the dynamics of the probability density of the vector
of phases, stochasticity arising only from the distribution of the natural frequencies.
They find, for fixed $N$, that there are stationary distributions of the vector of phases
clustered around a particular vector. We speculate that results of this type might
also hold for our very different model, \eqref{GMW1} and \eqref{GMW2}, which is based on quasi-cycles.

There is a significant mathematical literature, not focused on neuron populations,
about Kuramoto-type synchronization in the presence of noise. Two such recent
papers are by Giacomin et al. \cite{Gia12} and by deVille \cite{deville12}. Giacomin et al. studied
long-term dynamics of infinitely-many identical, noisy, phase oscillators using the
Kolmogorov equation of the stochastic system. When coupling strength is below
a critical value there is an incoherent state, e.g., phases are uniformly distributed.
When coupling strength is above this critical value the global attractor includes
a partially synchronized state. One would conjecture that similar results would
hold for infinitely-many interacting quasi-cycles. Even the careful statement of an analogue to the
Giacomin et al. results, however, is beyond the scope of this paper.

Finally, deVille~\cite{deville12} studied the long-term behaviour of a Kuramoto system of oscillators
perturbed by noise using Friedlin-Wentzell theory. That work concentrates on, but is
not limited to, nearest-neighbour interactions. This viewpoint is very different from ours, since we begin with a noisy system necessary to produce quasi-cycles, which are the objects we then regard as coupled. Even if we view our stochastic dynamic system, \eqref{GMW1} and \eqref{GMW2}, apart from its development here, we see that it cannot be regarded as a deterministic system perturbed by noise.
\section*{Competing interests}
The authors declare that they have no competing interests.
\section*{Author's contributions}
All authors contributed to the conceptualization and writing of the paper. The numerical simulations were
accomplished by LMW and MDM.
\section*{Acknowledgements}
Mark D. McDonnell is supported by the Australian Research Council under ARC grant DP1093425 (including an Australian Research Fellowship). Lawrence M. Ward was supported by Discovery Grant A9958 from NSERC of Canada.
  


\bibliographystyle{apsrev4-1}
\bibliography{gammabiblio.bib}

\end{document}